\documentclass{article}

\usepackage[cmex10]{amsmath}
\usepackage{amssymb}
\usepackage{graphicx}
\usepackage{cite}
\usepackage{url}
\usepackage{boxedminipage}
\usepackage{eurosym}
\usepackage{array}
\usepackage{color}
\usepackage{a4wide}

\newtheorem{rem}{Remark}
\newtheorem{lemma}{Lemma}

\def\bee{\begin{eqnarray*}}
\def\eee{\end{eqnarray*}}
\def\bmx{\begin{matrix}}
\def\emx{\end{matrix}}
\begin{document}

\title{Alleviating a form of electric vehicle range anxiety through On-Demand vehicle access}

\author{Christopher~King, %
     Wynita~Griggs, %
     Fabian~Wirth, %
     Karl~Quinn 
     and~Robert~Shorten
\thanks{This work was in part supported by Science Foundation Ireland: 11/PI/1177.}%
\thanks{C. King is with Northeastern University, Department of Mathematics.}
\thanks{W. Griggs is with the Hamilton Institute, National University of Ireland Maynooth, Maynooth, Co. Kildare, Ireland.}
\thanks{K. Quinn, F. Wirth and R. Shorten are with IBM Research Ireland, Dublin, Ireland.}}

\maketitle

\begin{abstract}
On-demand vehicle access is a method that can be used to reduce types of range anxiety problems related to planned travel for electric vehicle owners. Using ideas from elementary queueing theory, basic QoS metrics are defined to dimension a shared fleet to ensure high levels of vehicle access.  Using mobility data from Ireland, it is argued that the potential cost of such a system is very low.
\end{abstract}


\section{Introduction}

Growing concerns over the limited supply of fossil-based fuels are motivating intense activity in the search for alternative road transportation propulsion systems.  Regulatory pressures to reduce urban pollution, $CO_2$ emissions and city noise have made electric vehicles (EVs) \cite{b0} and plug-in hybrid vehicles (PHEVs) \cite{b1} a very attractive choice as the alternative to the internal combustion engine (ICE) \cite{rangeAnxiety}. Electric vehicles in particular, which have zero tailpipe emissions of pollutants,  are seen as a very useful tool in reducing urban pollution due to the fact that they deliver energy in a clean form in our cities,  as well as reducing the carbon footprint of road transportation generally when combined with clean energy generation.

While the environmental and societal benefits of zero-emission vehicles are evident, their adoption by users
has been extremely disappointing.  According to recent reports \cite{b4}, even in Europe, where the green agenda is well received,
fewer than 12000 electric vehicles were sold  in the first half of 2012 (of which only 1000 of these were sold in the UK). This number represents less than $0.15\%$ of total new car sales. These figures are in spite of the fact that many European governments have offered incentives for the purchase of electric vehicles in the form of subsidies and have also invested in enabling infrastructure. Some of the factors hindering the widespread adoption by the general public of electric vehicles can be summarised briefly as follows:

\begin{itemize}
\item[(i)] {\bf Battery related issues} -- Plug-in vehicles tend to be very expensive, even when subsidised. A major factor in the cost of such vehicles is the cost of the battery \cite{battery}. Lithium-based batteries are expensive, and while costs are forecast to reduce dramatically over the next few years \cite{b7,b8}, this is currently an important aspect in understanding the sales of electric vehicles. In response to this, Renault, and other companies, are proposing to lease batteries to the customer to offset some of the battery related costs. However, even if such initiatives are successful, there are other battery related questions that may hinder the adoption of electric vehicles. A further point concerns whether  enough lithium can be sourced to build batteries to construct enough vehicles to replace the existing passenger vehicle fleet. Are we simply substituting one rare resource (oil) with another (lithium)? Also, the transportation of batteries is not trivial and necessitates special precautions  \cite{shipping,b10}. Finally, most reasonably sized batteries are not capable of realising the range enjoyed by conventional ICE based vehicles.  While this latter issue is the subject of much research, battery size and performance currently represents one of the major determinants in the design of electric vehicles today \cite{ibm1,battery}.


\item[(ii)] {\bf Electromagnetic emissions} -- A recent issue regarding electric vehicles concerns electromagnetic emissions. While there is no evidence that EM radiation from EVs is dangerous, this issue is a focus point for regulatory authorities (see EU Green Car Programme) and has been raised by several research agencies \cite{b9}.


\item[(iii)] {\bf Long charging times} -- Charging times for electric vehicles
are known to be long \cite{b11}.  An often cited fact by advocates of electric vehicles in response to this is  that fast charging algorithms can service average vehicles in about 30 minutes \cite{b12,neiss}. Such time-scales may be just about acceptable to a normal car owner. However, in the presence of queuing, 30 minutes can rapidly become several hours, and push such fast charging stations into the realm of ``not acceptable.'' Thus, it is likely that overnight or workplace charging will be the principal method of vehicle charging for the foreseeable future. An associated issue in large cities concerns the availability of charging points. This is especially an issue in cities with large apartment block type dwellings.


\item[(iv)] {\bf Vehicle size} -- Electric vehicles tend to be designed small with limited luggage space to reduce energy consumption. This is a significant problem for most potential purchasers of vehicles who on occasion would like to transport significant loads using their vehicles.


\item[(v)] {\bf Range anxiety} -- One of the most pressing issues in the deployment of EVs concerns the issue of range anxiety \cite{rangeAnxiety}. Maximum ranges (in favourable conditions) of less than 150km are not unusual for electric vehicles, and this reduces significantly when air-conditioning or heating is switched on \cite{drivingRange}.  The issue of limited range also exacerbates other issues. For example, the cost of searching for a parking space at the end of a journey is much higher than for a conventional vehicle (because the EV's range is so low and therefore energy should not be wasted searching for a parking spot). Research is ongoing to address {range} issues, with much of the current work focussing on new battery types, optimal vehicle charging, vehicle routing, and in-vehicle energy management systems with a view to minimising wastage of energy and thereby increasing vehicle range \cite{shorten}. \protect\footnote{Much of this work mimics observed changes in driver behaviour when faced with the need to increase range. For example, in \cite{woodjack}, behavioural adaptations (in response to limited available energy) observed among participants of a study group, who were leased a battery EV for a year, were described.  Some of these behavioural adaptations included turning off the air conditioning or heater and  driving more slowly, as well as swapping vehicles with other users.}
\end{itemize}

While all of the above issues are important, as well as some others that we have not mentioned,  in this paper, we focus on the pressing issue of  {\bf range anxiety} (see below), and to a lesser extent, {\bf vehicle size} by suggesting a flexible vehicle access model to alleviate both of these issues.\newline

{\bf Principal contribution of paper :} Specifically, a solution is proposed to some of these problems based on a car sharing.  This idea is an embodiment of flexible vehicle access that was first suggested in \cite{ibmref} and further developed and analyzed in \cite{yale}. Indeed, the timeliness of the idea is evidenced by the fact that since this paper was first submitted for publication and presented \protect\footnote{The Sixteenth Yale Workshop on Adaptive and Learning Systems, Yale University, New Haven, Connecticut, USA, June, 2013.} \cite{ibmref, yale}, several automotive manufacturers, including Fiat, and Volkswagen, have announced plans for a form of car-sharing with similarities in both goals and implementation, to that described here. The principal difference to these embodiments and that proposed by IBM is that we advocate and describe an {\bf On-Demand model}, whereas car manufacturers are suggesting a model that guarantees access for a fixed, but small, number of days.  Our contribution in this paper is twofold. We develop tools to design an on-demand service using ideas from queueing theory and using predictive analytics. We then demonstrate,  in the context of real Irish mobility patterns, that such an On-Demand service is economical, both in terms of the number of 
ICE vehicles needed, and in terms of real additional cost to vehicle manufacturers.

We point out that in our queueing model it is assumed that extensive trips
are planned and the on-demand scheme responds to demand announced on the
previous day. There are, of course, other possible on-demand
scenarios. The assumptions described below pose no fundamental
restrictions in this respect and the methods presented here can be
extended.

\section{On types of Range Anxiety}
Our basic assumption is that the range anxiety problem, along with vehicle cost, are the biggest barriers
to the purchase of electric vehicles. While vehicle cost is principally an issue of manufacturing scale and is likely to
come down, range anxiety is an issue that is unlikely to be solved in a cost effective manner in the near future by developments
in the vehicle design process. For the purpose of
this paper we consider the term range anxiety to mean the {\em angst} of a vehicle owner that he/she will
not have enough range to reach his/her destination without the need for recharging. Roughly speaking, one may consider two types
of problems associated with range anxiety.
\begin{itemize}
\item[(i)] The first problem is associated with an inadequate battery level to complete a trip while driving to a destination.

\item[(ii)] The second problem is associated with electric vehicles being unsuited to the trip distribution demanded by the user.
For example, a long trip, or vacation (or indeed a trip where a large luggage load is required), are all problematic for a
typical electric vehicle.
\end{itemize}

In many cities the first of these problems can to a large extent be avoided with adequate trip planning. For example, in cities
where single dwelling households with garages/driveways are common, it is possible to charge vehicles adequately overnight
to have a full battery charge the next morning. As we shall see later, this
is, in most cases, adequate for the majority of trips. Since, most UK/Irish cities are of this type,
we shall assume that overnight charging is possible and focus on the second of the above range anxiety issues.

It is worth noting that the issue of range anxiety has been the subject of the attention of policy bodies and car manufacturers over the past number of years.
Roughly speaking, research attention has focussed on three areas: (i) better batteries; (ii) optimizing energy management
in the vehicles; and (iii) novel energy delivery strategies for electric vehicles. Big efforts on designing better batteries have been made worldwide  \cite{ibm1}.
To optimize energy management, manufacturers have looked both within the electric vehicle through the
management of the vehicle sub-systems, and outside of the vehicle
through the use of energy-aware routing and the use of special lanes (see www.green-cars-initiative.eu/projects/eco-fev). Finally, several methods of delivering energy to the vehicle
have been suggested. These include battery swapping, as was advocated by BetterPlace, under-road induction, and fast charging outlets. These activities
have addressed for the most part Item (i) above; and largely ignore the inconvenience associated with Item (ii).

\section{Problem Statement}\label{s1}

Our car sharing concept closely follows the flexible access suggestion in \cite{ibmref} in the following manner. \medskip

\noindent\begin{boxedminipage}[H]{1\linewidth}
When an electric vehicle is purchased, the new EV owner also automatically becomes a member of a car sharing scheme, where a shared vehicle may be borrowed from a common pool on a 24hr basis.  The shared vehicles are large ICE-based vehicles suitable for long range travel and with large goods transportation capacity.
\end{boxedminipage}
\medskip
\begin{rem}
We suggest free membership of the scheme, but a pricing model could be implemented to regulate demand on weekends, public holidays, or other occasions when synchronised (correlated) demand is likely to emerge, or to regulate emissions.  Further, if the shared ICE-based vehicles are chosen to be sufficiently high-end, then a further incentive for consumers to purchase electric vehicles is provided.
\end{rem}



\begin{figure}[h]
\begin{center}
\includegraphics[scale=0.36]{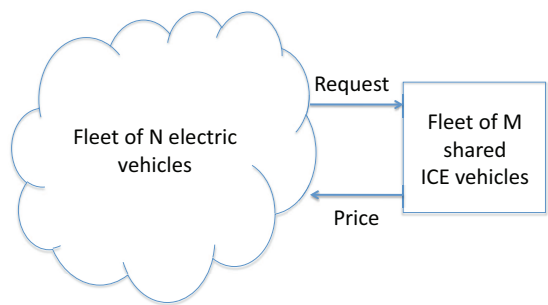}
\caption{Car sharing concept.}\label{fig:concept}
\end{center}
\end{figure}

A number of issues need to be resolved before any such system could be deployed.  These issues reduce to the marginal cost of the system. More specifically, we wish to determine if such a sharing concept could be deployed giving reasonable Quality of Service (QoS) to the electric vehicle owner, without significantly increasing the cost of each vehicle. Referring to Figure \ref{fig:concept}, this amounts to asking whether a reasonable QoS can be delivered when $M$, the number of shared ICE-based vehicles, is significantly less than $N$, the number of purchased EVs. To answer this question, we consider two scenarios.


\begin{itemize}
\item[1] {\bf  Spontaneous journeys :} If we want to achieve a low probability for the event that an individual
requesting a vehicle is not allocated one, how should the proportion $\frac{M}{N}$ be chosen?

\item[2] {\bf Planned journeys :}  For a fixed $M$ and $N$, how many days in advance does a user have to make a reservation so that the probability that the request is declined is lower than some very small constant?
\end{itemize}
In solving the above problems we shall make the following assumptions.\newline
{\em I. Membership assumptions}
\begin{itemize}
\item[I.1] A shared vehicle is borrowed for a single ``day" period and returned to the pool of shared vehicles to be shared again the next day.
In other words, at the beginning of each ``day", $M$ cars are available for sharing.
\item[I.2]  Customers collect their vehicles from a number of sharing-sites. They either park their own EV at the site when collecting the ICE vehicle,
or travel there using some other form of transport.
\item[I.3] Members of our car sharing scheme do not have access to an ICE based vehicle other than through the car-sharing scheme.
\item[I.4] Members of the car-sharing community are willing to accept the same $QoS$ metric for vehicle access.
\end{itemize}
{\em II. Demographic assumptions}
\begin{itemize}
\item[II.1] Long journeys in private cars are rare, meaning that the range of an electric vehicle, even under worst-case conditions (eg: air conditioning use, traffic congestion, bad weather), should be sufficient for most journeys.
\item[II.2] Most urban dwellings are houses rather than apartments, meaning that there are no structural impediments to overnight charging, and that a full overnight charge should be sufficient to satisfy the needs of most daily mobility patterns.
\end{itemize}
The validity of these latter assumptions are the subject of the next section.
\subsection{Data analysis and plausibility of assumptions}
To make a case for the plausibility of {Assumptions II.1 and II.2} we examined publicly available data on contemporary Irish mobility patterns.

\begin{figure}[h]
\centering
\includegraphics[scale=0.38]{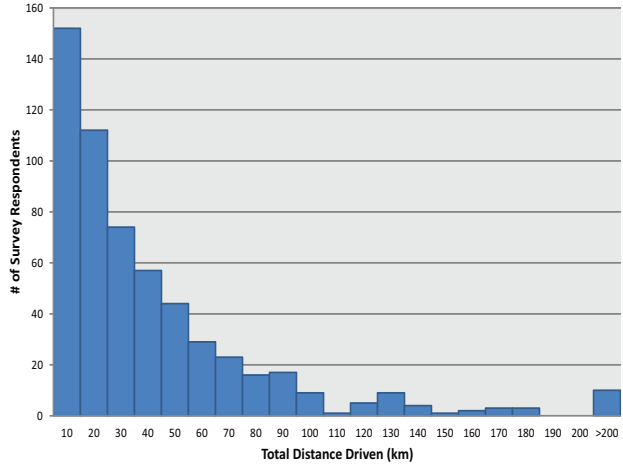}
\caption{Number of survey respondents reporting about their travel for the 24h ``Monday'' period, versus total distances they drove over that period.}\label{f1}
\end{figure}

\begin{figure}[h!]
\centering
\includegraphics[scale=0.38]{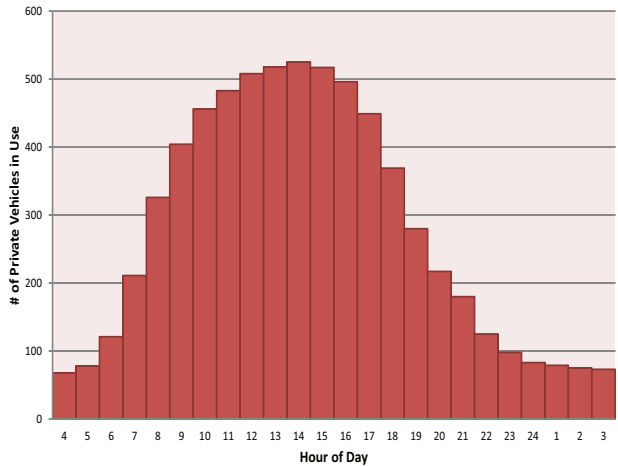}
\caption{Number of respondents (who drove $>$75km total daily distance) using their vehicles, versus hour of the day.}\label{f2}
\end{figure}

\begin{figure}[h]
\centering
\includegraphics[scale=0.38]{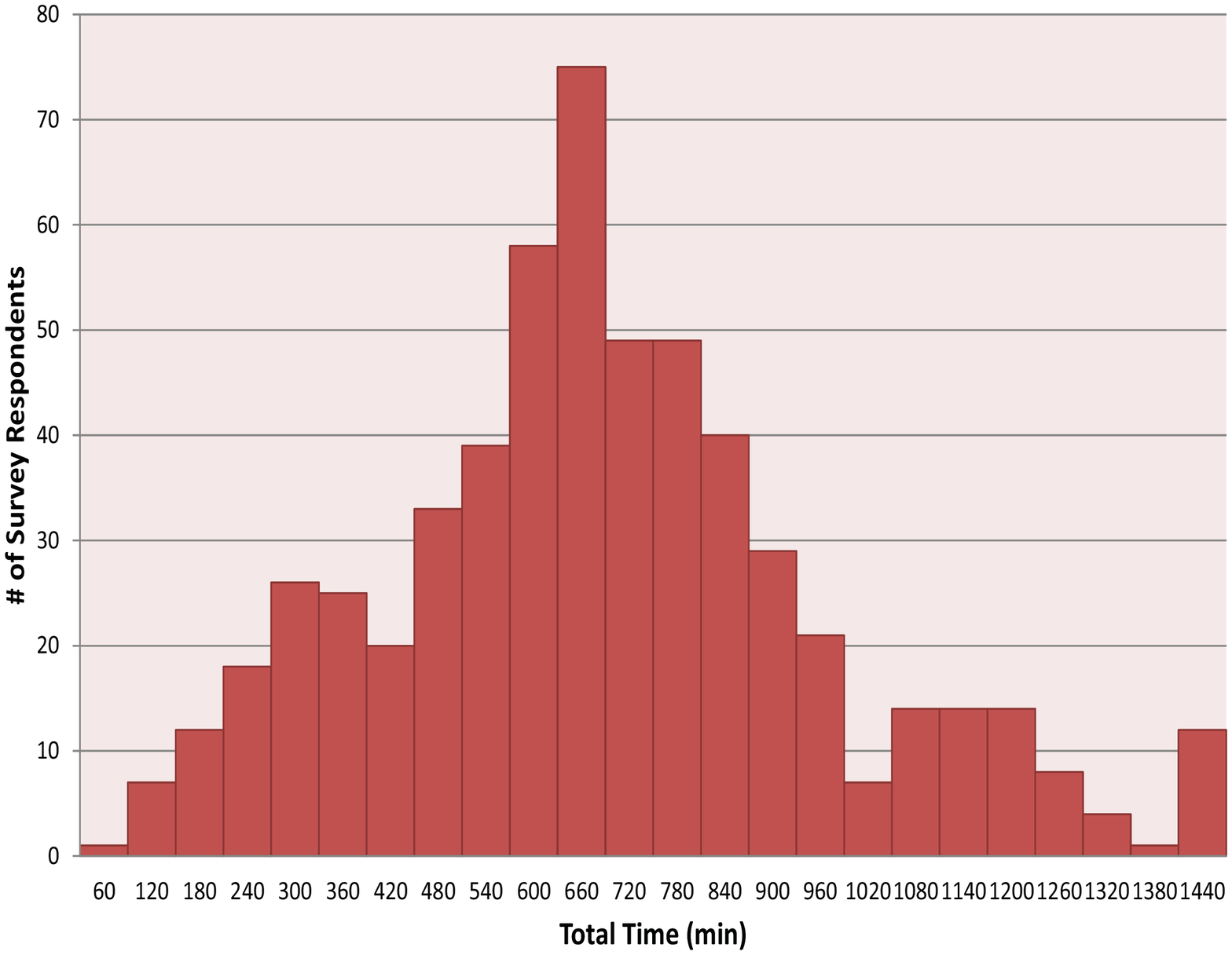}
\caption{Number of respondents (who drove $>$75km total daily distance), versus total time that their vehicle was in use for over their 24h survey period.}\label{f3}
\end{figure}

Data for the creation of Table \ref{tab1} and Figures \ref{f1}, \ref{f2} and \ref{f3} were obtained from the 2009 Irish National Travel Survey (NTS) Microdata File, Central Statistics Office, \copyright Government of Ireland \cite{CSOreport,CSO}. In the NTS, respondents were asked to provide details about their travel for a given (randomly selected) 24h period, which roughly corresponded to a day of the week.

Table \ref{tab1} shows the percentages of people who drove private vehicles (over the 24h period that they were queried about) for cumulative daily distances of greater than 50km, 75km and 100km. Figure \ref{f1} relates to those people who were questioned about their travel over the 24h ``Monday'' period (ie: row two of Table \ref{tab1}), and depicts number of people versus total distances they drove in private vehicles over that 24h Monday period. Figure \ref{f1} illustrates a trend observed in the percentages in Table \ref{tab1}; namely, that longer cumulative journeys over the course of a day were rare\protect\footnote{Graphs of the nature of Figure \ref{f1}, but concerning travel over the other days of the week, were similar in shape to Figure \ref{f1}, and have thus been omitted.}.  For respondents who drove cumulative distances greater than 75km over a 24h period (see the third column of Table \ref{tab1}), Figure \ref{f2} illustrates those hours of a 24h period over which respondents had their vehicles in use (many vehicles were in use roughly between 8am and 6pm), and Figure \ref{f3} depicts number of respondents versus total time (out of a 24h period) their vehicle was in use.

\newcolumntype{s}{>{\centering\arraybackslash}m{0.8cm}} 
\newcolumntype{u}{>{\raggedright\arraybackslash}m{2.2cm}} 

\begin{table}[h]
\caption{Percentages of people who drove cumulative distances of greater than 50km, 75km and 100km over a 24h period.}
\begin{center}
\begin{tabular}{u|sss}\hline\hline
\smallskip Sample Population & \smallskip 50km & 75km & 100km \\ \hline
\smallskip Monday & \smallskip 23\% & \smallskip 12\% & \smallskip 7\%\\
\smallskip Tuesday & \smallskip 23\% & \smallskip 14\% & \smallskip 8\%\\
\smallskip Wednesday & \smallskip 23\% & \smallskip 14\% & \smallskip 7\%\\
\smallskip Thursday & \smallskip 26\% & \smallskip 18\% & \smallskip 11\%\\
\smallskip Friday & \smallskip 26\% & \smallskip 17\% & \smallskip 9\%\\
\smallskip Saturday & \smallskip 24\% & \smallskip 15\% & \smallskip 9\%\\
\smallskip Sunday & \smallskip 24\% & \smallskip 17\% & \smallskip 11\% \\ \hline\hline
\end{tabular}
\end{center}\label{tab1}
\end{table}

\subsection{Comments on NTS data-set}
The NTS data set allows us to make some very general statements which support the assumptions listed above. As might be anticipated, longer journeys are not as frequent as shorter ones, and vehicles do indeed tend to be in-use for a single day for each user. Thus, the data set appears to be indeed consistent with the assumptions II.- II.2 above. However, the data-set is subject to certain limitations due to its size and the manner in which it was collected, and it is very important to comment on these before proceeding.

\begin{itemize}

\item First, the data-set does not take into account correlated driving behaviors between days. However, we assume that drivers making regular long-distance trips would not purchase EV's. Consequently, the tail of the distributions following from the data-sets should be representative of the qualitative behavior of drivers making regular short trips, but who infrequently need larger vehicles or vehicles for larger trips. In what follows we shall assume that the typical EV owner follows a distribution given by these data sets - although clearly from the above discussion this is a simplification as driving behavior varies across EV owners.

\item Second, it is possible that the car-sharing scheme may affect the driving patters of scheme members (EV owners). For example, given ICE vehicle access, long-distance driving may become more attractive. In the analysis to follow (in the next section), this feedback mechanism is also ignored to keep matters simple.

\item The data also indicates that vehicles are ``in-use'' for the entire working day. We assume that vehicles are charged overnight to full capacity and that it is not necessary for a vehicle to be again charged during the day. As we have already mentioned, this latter assumption is consistent with cites such as Dublin where garage space or a driveway is the norm, and is consistent with charging-times for standard home charge points.  The interested reader is referred to  \cite{ref3} for further information about driver charging behavior.
\end{itemize}

Given the above discussion it is clear that the presented data provides a plausibility argument in support of assumptions II.1 and II.2 subject to the aforementioned limitations of the data-set. For a further and more detailed discussion of the limits of using such data sets the interested reader is referred to \cite{ref4}.

\section{Mathematical Models}\label{s2}
Using elementary probability and queueing theory methods, we now pose solutions to the Quality of Service (QoS) problems presented in Section \ref{s1}.  Consider a population of $N$ electric vehicle owners (ie: $N$ ``users'') who occasionally require access to an ICE-based vehicle (ICE or ICEV) for a non-standard trip (either a long-range trip or a trip where large load carrying capacity is required).  We assume that a user will keep the ICEV for a full day, based on the driver behaviours described in the previous section.  Thus, each day is characterised by the number of users who require an ICEV on that day.  There is a fleet of $M$ ICEVs available to satisfy this need.  The main question is then to determine the relation between $M$ and $N$.  This will be determined by requiring some QoS conditions to be met.  For example, a QoS condition might be a guarantee on the probability of finding an ICEV available.

An important assumption in what follows is that all users are willing to accept the same $QoS$ metric for vehicle access. This is defined by a probability of a scheme member  not being able to access a vehicle when required. Similar assumptions are standard in the networking community, and are based on the assumption that a member of the car-sharing scheme will not join the scheme unless the given QoS metric is acceptable.  Clearly, different levels of service can be provided if groups of users are willing to pay for a higher level of service. However, to keep matters simple here we assume that drivers who sign up for the scheme are prepared to be served using the same QoS metric.

\subsection{Model 1 -- Binomial Distribution}

In the simplest model, each user independently requests an ICEV each day
with probability $p$. Based on the data available in Table \ref{tab1} we may
estimate this probability as $p=.0886$, assuming that ICEVs are used for
trips over 100km. The assumption that cars are required for a whole day
are justified by the data in Figures~2-4, which show that the
overwhelming majority of cumulative trips over 75 km require a period
of 5 hours or more. With the available data we are not able to estimate
the proportion of drivers that make long trips on a regular basis, who
would not participate in the scheme. On the other hand it may be expected that some ICEVs are requested even though the
actual trip length turns out to be less than 100km. In the following we will use
$p=.1$ to illustrate the applicability of the proposed scheme.

Thus, the number of requests $X$ each day is a binomial random variable:
\bee
X \sim {\rm Bin}(N, p).
\eee
The mean number of requests per day is $N p$, and the standard deviation is $\sqrt{N p(1 - p)}$. In principle, the number of requests may be anything from $0$ to $N$, but for large $N$ it is very unlikely that $X$ will deviate from the mean by more than a few standard deviations.


This QoS condition can be quantified as follows. For each $M \le N$, define
\bee
Q(M) = P(X > M).
\eee
Then the QoS condition could be to find the smallest $M$ such that $Q(M) < \epsilon$ for some specified $\epsilon$.
For any given $N$ and $p$ this can be calculated explicitly using the formula
\bee
Q(M) = \sum_{k=M+1}^N \left(\bmx N \cr k \emx\right) \, p^k \, (1-p)^{N-k}.
\eee
However it is more useful to get an approximate formula from which the scaling relation can be read off.
For $N$ large enough we can use the normal approximation for the binomial, which says that
\bee
Z = \frac{X - N p}{\sqrt{N p (1-p)}}
\eee
is approximately a standard normal random variable.
There is a standard rule of thumb regarding applicability of the normal approximation for the binomial, namely that
\bee
N \ge 9 \, \max \left(\frac{p}{1-p}, \frac{1-p}{p} \right).
\eee


Using the normal approximation, we have
\begin{equation*}
Q(M) \simeq  \frac{1}{\sqrt{2 \pi}} \, \int_{r}^{\infty} e^{- \frac{1}{2} X^2} \, d X,
\qquad r = \frac{M - N p}{\sqrt{N p (1-p)}}.
\end{equation*}
This readily yields estimates for $M$ in order to satisfy a desired QoS condition. For example, in order to satisfy the
QoS condition
\bee
Q(M) < 0.05
\eee
meaning a less than $5 \%$ chance of not finding an ICEV available, it is sufficient to take
\bee
r \ge 1.65 \Longleftrightarrow M \ge N p + 1.65 \sqrt{N p (1-p)}
\eee
For example, using the values $N=1000$ users and $p=0.1$ for the probability of a user requesting an ICEV, this
provides a value $M \ge 116$.

\subsection{Model 2 -- A Queueing Model}

When the number of requests exceeds the number of available ICEVs, a queue forms and users must wait one or several
days until a vehicle becomes available. It is desirable to keep the probability of long delays small, and this can be achieved by
appropriate scaling of $M$ with $N$. This is the subject of this section.


Let $X_n$ be the number of outstanding requests at the end of the $n^{th}$ day,
and let $A_n$ be the number of new requests that arrive during the $n^{th}$ day. Since the number of vehicles is $M$,
the relation between these variables on successive days is
\bee
X_{n+1} = \max \left\{ 0, X_n - M  \right\} +  A_{n+1}.
\eee
That is, the queue length is reduced by $M$ at the start of each day (but not reduced below zero), and is then increased
during the day by the number of new requests.


The QoS condition is to ensure that $X_n$ is unlikely to be large, implying that users are unlikely to have to wait a long time before
being assigned an ICEV. Since $M$ users can be serviced each day, in the worst-case a user must wait $\lfloor X_n/M \rfloor$ extra
days until service.
We consider the QoS condition which guarantees that the probability that any user needs to wait $k$ extra days or more
is less than $\epsilon$, that is
\bee
P(X_n > k M) < \epsilon.
\eee
By choosing $M$ sufficiently large we can guarantee that this probability is small.

\medskip

\noindent\begin{boxedminipage}[htb]{1\linewidth}
\begin{lemma}\label{lem:prob-bound}
Define
\bee
\mu = M - N p, \qquad \sigma^2 = N p (1 - p), \qquad \alpha = \frac{\mu}{\sigma^2}.
\eee
Then for all $k \ge 1$,
\bee
P(X_n > k M) \le \frac{1}{2} \, e^{-(k-1) M \alpha} \, \left(e^{\mu \alpha/2} - 1\right)^{-1}.
\eee
\end{lemma}
\end{boxedminipage}

\medskip

\begin{figure}
\begin{center}
\includegraphics[scale=0.45]{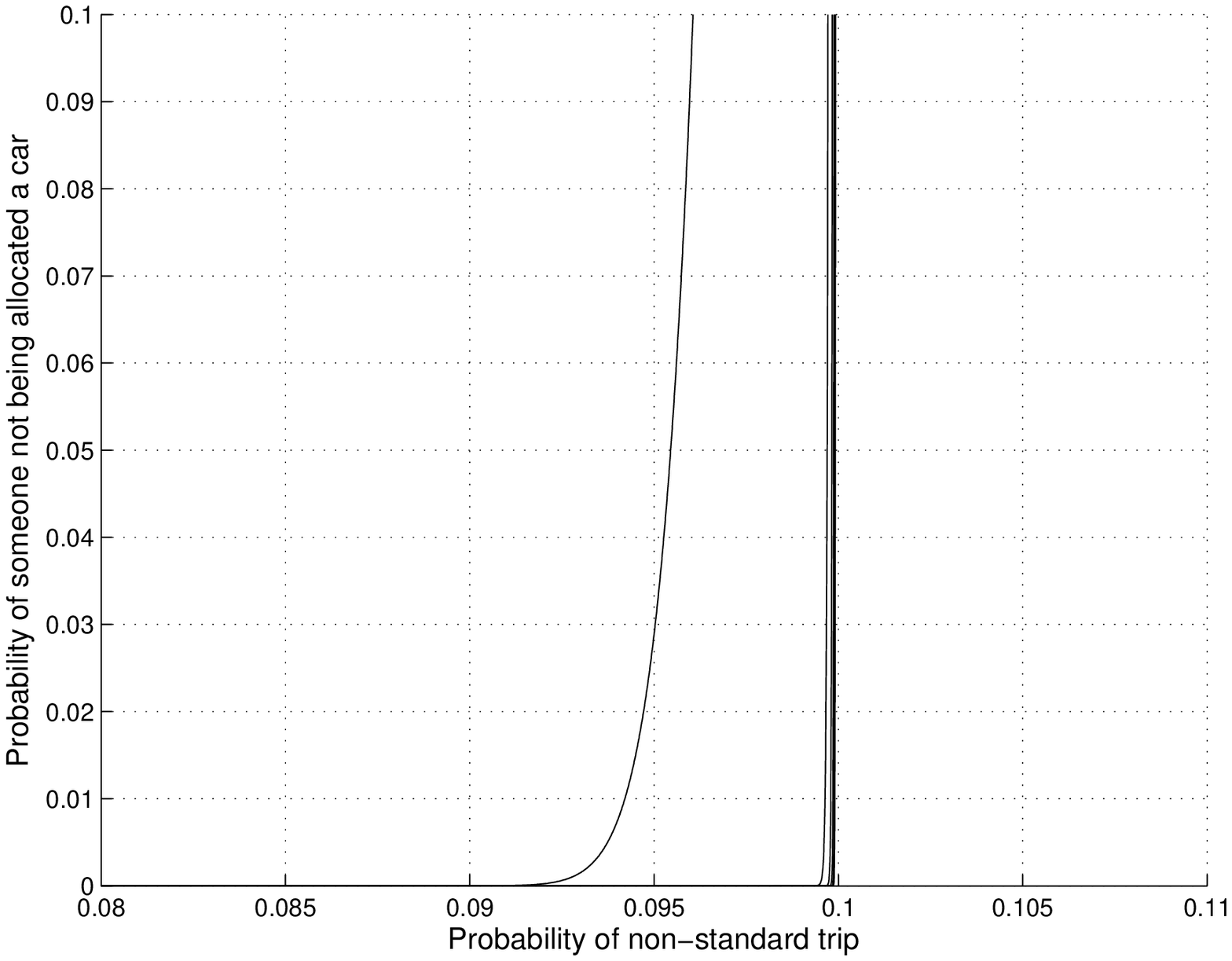}
\caption{Probability (bound) of not finding a car with N = 20000 and M = 2000 (10\%).  $k =1$ is the left most curve; $k=6$ is the right most curve.}\label{fig:model2}
\end{center}
\end{figure}

\begin{figure}
\begin{center}
\includegraphics[scale=0.45]{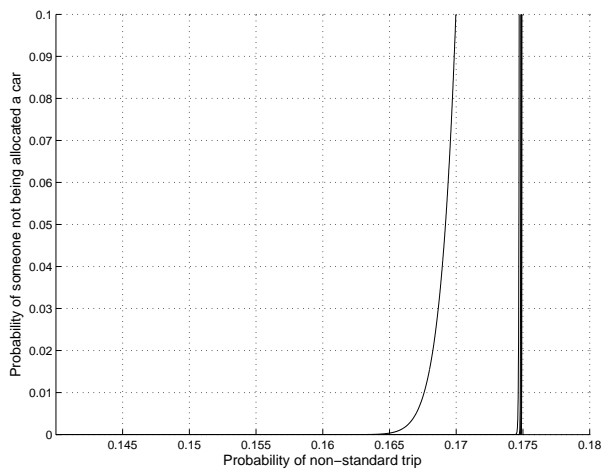}
\caption{Probability (bound) of not finding a car with N = 20000 and M = 3500 (17\%). $k =1$ is the left most curve; $k=6$ is the right most curve.}\label{fig:model3}
\end{center}
\end{figure}

Using the bound in Lemma \ref{lem:prob-bound} we find a sufficient condition to guarantee the QoS bound, namely
\bee
\frac{1}{2} \, e^{-(k-1) M \alpha} \, \left(e^{\mu \alpha/2} - 1\right)^{-1} < \epsilon.
\eee
For a given $k$ and $\epsilon$ we may use this to find a value for $M$ needed to meet the QoS condition.
For example, using the same values as above $N=1000$, $p=0.1$, $\epsilon = 0.05$, and taking $k=3$
(meaning that the probability that there is a customer who waits more than 4 days is less than $5\%$), we find that the inequality is
satisfied whenever $M \ge 103$. Taking $k=2$ we find $M \ge 105$, and with $k=1$ we find $M \ge 121$.
For a fleet size of $20000$ vehicles with $M=2000$ $(10\%)$ and $M=3500$ $(17\%)$, the behaviour of the bound is depicted in Figures
\ref{fig:model2} and \ref{fig:model3}, respectively, for various estimates of the probability of a long distance trip.
As can be seen, the bound tends rapidly to zero, indicating that the probability of waiting for a shared vehicle longer than one or two days vanishes rapidly.


{\it Proof of Lemma \ref{lem:prob-bound}:} Define
\bee
W_n = \max \left\{0, X_n - M \right\}, \qquad B_n = A_n - M.
\eee
Then the recursion relation becomes
\bee
W_{n+1} = \max \left\{0, W_n + B_n  \right\}.
\eee
As is well-known \cite{klienrock}, with $W_0=0$ the solutions of this recursion relation have the same (marginal)
distributions as the sequence $W_0', W_1', \dots$ where $W_0'=0$ and
\bee
W_n' = \max \left\{0, B_1, B_1 + B_2, \cdots \right\}.
\eee
So define
\bee
C_j = \sum_{i=1}^j B_i.
\eee
Then we have for $k \ge 1$
\begin{align*}
P(X_n > k M) &= P(X_n - M > (k-1) M ) \\
&= P(W_n > (k-1) M) \\
&= P ( \max \left\{0,C_1,C_2,\dots,C_n \right\} > (k-1) M) \\
&= P \left( \bigcup_{j=1}^n \left\{C_j > (k-1) M \right\} \right)\\
& \le \sum_{j=1}^n P(C_j > (k-1) M).
\end{align*}
For large $N$ we can use the normal approximation to estimate the right side of this inequality.
The arrival numbers $A_n$ are independent and  binomial with mean $N p$ and variance $N p(1-p)$, so in this
approximation the variables $C_j$ are also normal, with mean $j (N p - M)$ and variance $j N p (1-p)$.
Thus we find
\begin{equation*}
P(C_j > (k-1) M) = P\left(Z > \frac{(k-1) M - j (N p - M)}{\sqrt{j N p (1-p)}} \right)
\end{equation*}
where $Z$ is a standard normal. We use the following bound valid for all $b \ge 0$
\bee
P(Z > b) \le \frac{1}{2} \, e^{- b^2/2}
\eee
to obtain
\bee
P(C_j > (k-1) M) \le \frac{1}{2} \, e^{ - ((k-1)M + j \mu)^2/2 j \sigma^2}.
\eee
We use the bound
\bee
\left(\frac{x}{\sqrt{j}} + \sqrt{j} y \right)^2 \ge 2 x y + j y^2
\eee
and apply it above with $x = (k-1)M $ and $y = \mu$.  This gives
\bee
P(X_n > k M) & \le & \sum_{j=1}^n \frac{1}{2} \, e^{- (k-1)M \mu/\sigma^2} \, e^{ - j \alpha \mu/2}.
\eee
The sum is geometric, and is bounded by
\bee
\sum_{j=1}^n e^{ - j \alpha \mu/2} < \sum_{j=1}^{\infty} e^{ - j \alpha \mu/2} = \frac{e^{ -  \alpha \mu/2}}{1 - e^{ - \alpha \mu/2}}
\eee
as required. $\blacksquare$

\section{Financial Calculations}\label{s3}

We now explore the financial model associated with this car sharing model. Initially, we will model our calculations based on a shared fleet comprised solely of a single vehicle type: the VW Golf class car. In this model we are essentially comparing the similarly sized VW Golf and the Nissan Leaf to focus solely on examining range anxiety
of the kind discussed above.

Subsequently, we will model our calculations based on a shared fleet that is constructed to reflect the broader needs of the EV owners; namely, sometimes a family car is needed, sometimes a large vehicle is needed for transporting goods, and sometimes a smaller car is required for short out of town trips by a single person/couple. In this subsequent model we focus on both range anxiety and a range of vehicle sizes. To do this we make some simplifying assumptions. Instead of modeling the demand for various vehicle types, which is the correct method of analysis, we shall weight various vehicle classes to construct a vehicle fleet. We assume that these weights reflect demand. This simplifying model is adopted for two reasons; to keep analysis simple, and to reflect the fact that data of the type needed to build a multi-dimensional queueing model is not available to us at this time and is the subject of current work in our lab.

In both model calculations we use Volkswagen (VW) vehicles. Note that VW is a relatively high-end marque, and the vehicle fleet could be constructed in a manner that is considerably cheaper. Table \ref{table:nonlin} illustrates the composition and pricing of our shared vehicle fleet (including the single vehicle VW Golf price). Vehicle purchase prices were sourced from the VW website (www.volkswagen.ie) on 23rd September, 2013. The weighted average vehicle purchase price given the fleet make-up is \euro 23,308.

\begin{table}[h]
\caption{Composition of the shared vehicle fleet.}
\begin{center}
\begin{tabular}{l c l c}\hline\hline
Purpose & $\%$ of Fleet & Vehicle Purchase & Price  \\
\hline
Singles/Couples &20\%& VW Polo &  \euro 14,195\\
Family without luggage & 20\%&VW Golf & \euro 19,995 \\
Family with luggage  & 50\% & VW Passat  &  \euro 27,165 \\
Transport & 10\% & VW Passat Estate  &  \euro 28,870 \\
\hline
\end{tabular}
\end{center}\label{table:nonlin}
\end{table}

We also assume a fleet of $N$ Nissan Leaf electric vehicles. These electric vehicles currently retail for \euro 25,990 each in Ireland (23/9/2013).
It is important to note that there are more expensive and less expensive ICEV and EV offerings.

Ahead of detailing the two model calculations we recall Table \ref{tab1}. Table \ref{tab1} indicates that the probability of cumulative journeys greater than 75 km is approximately 0.15 (average), and of 100 km is 0.09 (average). If we conservatively assume that the daily range of a fully charged electric vehicle is 75 km, then it follows from Figure 6 that most customers will be allocated an ICE-based vehicle within 3 days for (M,N) = (3500, 20000); namely, if M/N = 0.17. If we assume that the daily range is 100 km, then it follows from Figure 5 that most customers will be allocated an ICE-based vehicle within 3 days for (M,N) = (2000, 20000); namely, if M/N = 0.10. In essence, in order to provide 20,000 EV cars with an ICEV car within 3 days then we need an additional 2,000 ICEV cars to cater for journeys of 100km and an additional 3,500 ICEV cars to cater for journeys of 75km.

\subsection{Range Anxiety Model (VW Golf vs. Nissan Leaf)}

We utilise a financial model for both the car sharing models, which is in line with current and well known car financing structures. We take an approach that amortizes the purchase of a new ICEV over a three year period with a 20\% straight line depreciation year on year. However, first year depreciation is at 40\% as suggested as standard by AA (http://www.theaa.com/motoring$\_$advice/car-buyers-guide/cbg$\_$depreciation.html). This gives three year depreciation costs for a VW Golf as listed in Table \ref{table3}.

\begin{table}[h]
\caption{Depreciation costs of the VG Golf fleet.}
\begin{center}
\begin{tabular}{l c c c c}\hline\hline
& Start Value & Year 1 & Year 2 & Year 3  \\
\hline
Car Value & \euro 19,995 & \euro 11,997 & \euro 9,598 & \euro 7,678 \\
YoY Depreciation (\%) & & 40\% & 20\% & 20\% \\
YoY Depreciation (\euro) & & \euro 7,798 & \euro 2,399& \euro 1,920\\
$\%$ Value vs. Start & & 60\% & 48\% & \textbf{38\%} \\
GMFV & & & & \textbf{40\%}\\
\hline
\end{tabular}
\end{center}\label{table3}
\end{table}
It is important to note that we also make the assumption that the ICEV is sold in year three for a value close to its depreciated value. Motor companies such as Hyundai, Ford, etc. provide Guaranteed Minimum Future Value (GMFV) mechanisms such that the depreciated value of the car is close to the GMFV of the car after a given time that is typically 3 years. The three year amortization of depreciation and GMFV allow us to financially model the depreciation costs of the ICEV within a 3 year bound.

Considering the revenues from 20,000 EV cars at a cost of \euro 25,990 is \euro 519.8 million we will now present the additional cost of this car sharing model using the 3 year financial approach outlined above.

In the case of 100km journeys we have described how 2,000 ICEV are required. As we are using VW Golf in this model then the 2,000 ICEVs would have once off costs of \euro40m, which is 7.7\% of the EV revenues. However, such a once off financing mechanism is not typical of financing car purchases in the car industry. Instead, we utilize the 3 year amortization and GMFV method outlined above so as the costs can be modeled more realistically. In Table \ref{table3} the sum of the single VW Golf depreciation values (YoY Depreciation \euro) over the three year period is \euro12,317. In year three we dispose of the VW Golf asset at, or near, the GMFV value of ~40\%. As such, the cost of the VW Golf over three years then becomes the depreciation costs of \euro 12,317. As we require 2,000 VW Golf cars for this car sharing model this then equates to a fleet cost of \euro 24.6m over three years, or an average of \euro 8.2m per annum. This represents an annual cost overhead of 1.5\% per annum for three years for 100km journeys serviced by a fleet of 2,000 VS Golf cars within 3 days. In the case of 75km journeys we require 3,500 VW Golf cars, which has a similarly calculated annual cost overhead of 2.7\%.

\subsection{Range Anxiety with a Range of Vehicle Sizes Model}

This model will utilise the same 3 year financial approach for a range of ICEV as used in the VW Golf model. However, the cost of the ICEV is the weighted average, \euro 23,308, of a range of VW ICEV types taken from Table \ref{table:nonlin}. In this scenario we are modeling both range anxiety and the ability to cater for ICEV usages for a family (i.e. VW Passat) with luggage and transport purposes (i.e. VW Passat Estate).

\begin{table}[h]
\caption{Depreciation costs of the weighted ICE fleet.}
\begin{center}
\begin{tabular}{l c c c c}\hline\hline
& Start Value & Year 1 & Year 2 & Year 3  \\
\hline
Car Value & \euro 23,380 & \euro 13,985 & \euro 11,188 & \euro 8,950 \\
YoY Depreciation (\%) & & 40\% & 20\% & 20\% \\
YoY Depreciation (\euro) & & \euro 9,323 & \euro 2,797& \euro 2,238\\
$\%$ Value vs. Start & & 60\% & 48\% & \textbf{38\%} \\
GMFV & & & & \textbf{40\%}\\
\hline
\end{tabular}
\end{center}\label{table4}
\end{table}

The revenues from 20,000 EV cars at a cost of \euro 25,990 remain unchanged at \euro 519.8 million.

As per the previous model 100km journeys require 2,000 ICEV. As we are using a weighted average of different ICEV class cars in this model then the 2,000 ICEVs would have once off costs of \euro 46.6m, which is 9\% of the EV revenues. In Table \ref{table4} the sum of the depreciation values (YoY Depreciation \euro) over the three year period is \euro 14,358. As we again require 2,000 ICEV then this equates to a fleet cost of \euro 28.7m over three years, or an average of \euro 9.6m per annum. This represents an annual cost overhead of 1.8\% per annum for three years for 100km journeys that are serviced by a range of ICEV within 3 days. In the case of 75km journeys we require 3,500 ICVEs, which has a similarly calculated annual cost overhead of 3.1\%.

\subsection{Financial Assumptions and Key Conclusions}

In both models we assume further incremental costs such as traveling to ICEV location, parking logistics, car cleaning, annual servicing, and all other operational costs are either mitigated via an implemented pricing model and/or absorbed by the fixed cost structure already with the dealership, rental company, and so on. We also assume that at the end of the three year period the EV car owner must stop using the service or replace that EV with a new EV to start another 3 year cycle, which maintains the calculations and logic of the approach over terms greater than an initial 3 years.

In the model where we used VW Golf to focus solely on alleviating range anxiety the overall percentage cost is between 4.5\% (100km, 3 years) and 8.1\% (75km, 3 years) of the 20,000 EV costs. In the weighted average model where we focused on alleviating range anxiety, and also offering a range of ICEV size options, then the overall percentage cost is between 5.4\% (100km, 3 years) and 9.3\% (75km, 3 years) of the 20,000 EV costs. Note, for further clarification the boundaries of our model ranges from 3.3\% overall costs (2,000 VW Polo, 100km) to 11.7\% overall costs (3,500 VW Passat Estate, 75km).

To place these figures into context, consider Table \ref{table5}. As can be seen, the cost of the car sharing scheme is considerably less than the average level of subsidy afforded to electric vehicles in major western countries; currently at approximately 23\% for a Nissan Leaf EV. Further incentives to encourage the uptake of electric vehicles might include a combination of car sharing/subsidies or replacing the subsidies with an increased level of car sharing. Whether car sharing can really encourage the uptake of such vehicles can only really be tested through implementation. However, we believe that we have demonstrated that a significant range related issue can be solved using this idea, and that this will ultimately affect market growth in a positive manner. 

\begin{table}[h]
\caption{Subsidies to EV purchase (direct and indirect) and cost of Nissan Leaf. Data sourced from Nissan: September 2013}
\begin{center}
\begin{tabular}{l c c c}\hline\hline
Country & Subsidy & Cost (Nissan Leaf) & Percentage \\
\hline
Ireland & \euro 5000 & \euro 25990 &  19\%\\
Belgium & \euro 9000 & \euro 29890&  30\%\\
France  &  \euro 7000 & \euro 30190 &  23\%\\
Portugal & \euro 5000 & \euro 31100 &  16\%\\
United Kingdom & \pounds 5000 & \pounds 20990 &   25\%\\
United States  & \$7500 & \$28800 &  26\%\\
\hline
\end{tabular}
\end{center}\label{table5}
\end{table}

To conclude it is worth recalling some of the assumptions underlying our analysis.

\begin{itemize}
\item Our analysis ignores synchronized demand (at weekends or during holiday periods). We argue that a pricing structure could be enforced to
breakup this demand and give a degree of QoS to the scheme members.

\item As we have mentioned, our combined fleet analysis is not based on demand distributions for each vehicle class. Rather, it is a simplified calculation
to give the reader an indicative picture of the cost of the scheme.

\item We have ignored the fact that the shared vehicles would be in continuous use thereby rendering their value lower than standard GMFV. However, GMFV mileage restrictions may be offset through gains made when bulk buying 2,000 or 3,500 ICEV, and/or through negotiating greater mileage restriction caps, and/or by ensuring a reduction in any additional charges beyond the restricted mileage, and so on. Additionally, other depreciation models are easily incorporated into the given framework. Note that depreciation of the entire cost of the shared vehicles over a 3 year period to zero value is also very cost effective when compared with the average cost of government subsidies (7.7\% for a VW Golf in the 100km model to 15.7\% for a weighted average ICEV in the 75km model, in comparison with the 23\% average subsidy). A fleet of 3,500 VW Passat Estates in the 75km model depreciated to zero is more efficient than subsidies (19.4\% compared to 23\%).

\end{itemize}
Notwithstanding the above facts, we have shown how any/all of the costs from the two model calculations are significantly less than the government subsidies costs. This gives governments and/or industry the opportunity to augment and/or replace subsidies with this alternative model to encourage the uptake of EV cars, make EV cars more approachable, and reduce pollution in short journeys.

\section{Final comments}
\begin{figure}
\begin{center}
\includegraphics[scale=0.45]{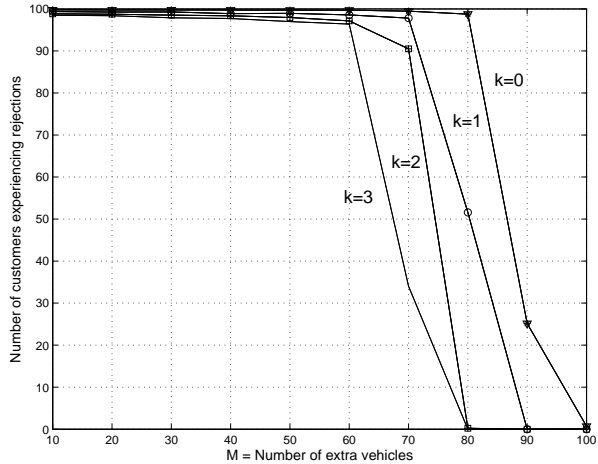}
\caption{A simulation to illustrate the behavior of the queueing model. For bookings greater than 3 days in advance, the number of unhappy customer goes to zero.}\label{fig:model4}
\end{center}
\end{figure}

{\bf Simulation :} To conclude our paper we now present a brief simulation to illustrate the
performance of our car sharing system. We simulated demand based on Figure 2, for 365 days,
and serviced this demand using our car-sharing model.  Users requested an ICE  (no
accommodation for type of ICE was implemented in our simulation) with a fixed probability of $p=0.1$ (to emulate infrequent need to take a long journey), independently of each other. The number of members
of the car-sharing program was $N = 1000$. The size of the ICE fleet was varied from $M = 10$ to $M=100$.
Figure \ref{fig:model4} depicts the percentage of customers who experienced
not being able to access a car within $k$ days over a 365 day period.
As can be seen, the number of unhappy customers (rejections)
goes rapidly to zero as $M$ approaches 100. This is entirely consistent with the conclusions
presented above and illustrates that the predictions of the mathematics are upper bounds on performance.\newline

{\bf Carbon savings :} Finally, to quantify, very roughly, the emissions savings over an ICE fleet, we consider 
the distance distribution depicted in Figure 2 under the following assumptions.  We assume 
that all journeys are made with a single vehicle type traveling at a constant speed.  Further,  we assume that long distances are those greater
than 100 $km$. Then, based on Figure 2, the CO$_2$ saving is roughly $80\%$ using the proposed car-sharing scheme. \newline

{\bf Collaborative consumption :} We believe that the techniques presented in the paper go far beyond electric vehicles. Many sharing schemes (sometimes called collaborative consumption schemes) reduce to asking how many devices/units are needed to serve a  demand, while at the same time delivering a statistical QoS to
individual users. We believe our ideas and models will help in the context of such applications.

\section{Conclusions and Future Work}\label{s4}

In this paper, a solution to a consumer range anxiety problem is studied using a car sharing idea. The cost of this scheme was shown to be low when compared with current levels of subsidies to electric vehicle manufacturers. This cost could be reduced further by exploiting car sharing on an hourly basis (further multiplexing) and using low-cost vehicles (as opposed to a premium marque).

\section*{Acknowledgment}

The authors thank Douglas Leith and Lisa Amini for their valuable comments and the excellent comments of the reviewers.
We would also like to thank Emanuele Crisostomi for noticing an error in Figures 5 and 6.

\end{document}